\documentclass[11pt]{article}
\usepackage{amsthm, amsfonts}
\usepackage{amscd,amsmath,amssymb}

\newtheorem{definition}{Definition}
\newtheorem{theorem}{Theorem}
\newtheorem{lemma}{Lemma}

\author{A.~M.~Vershik\thanks{%
St.~Petersburg Department of Steklov Institute
of Mathematics. Fontanka 27, 191023 St.~Petersburg, Russia.
E-mail: {\tt vershik@pdmi.ras.ru}.
Partially supported by RFBR under the project 05-01-00899
and INTAS under the project 03-51-5018.}}
\title{What does a generic Markov operator look like?}
\date{}

\begin{document}
\maketitle

\rightline {\it To the memory of O.~A.~Ladyzhenskaya}

\begin{abstract}

We consider generic (i.e., forming an everywhere dense massive subset)
classes of Markov operators in the space
$L^2(X,\mu)$ with a finite continuous measure. Since there is a canonical
correspondence that associates with
each Markov operator
 a multivalued measure-preserving transformation
(i.e., a polymorphism), as well as a stationary Markov chain,
we can also speak about generic polymorphisms
and generic Markov chains. It was not known
not only that the suggested properties are generic,
but even whether  there exist
Markov operators that have simultaneously all or a part of them.
The most important result is that the class of
totally nondeterministic nonmixing operators is generic.
We pose a number of problems and express the hope  that
generic Markov operators will find applications in various fields, including
statistical hydrodynamics.
\end{abstract}

I was lucky to be friends with O.A., especially in the 70s,
and some time I will write about this. In the late 60s and 70s, she
was very interested in dynamical systems, and this was an
additional motive for our contacts.
Working on the Hopf
equation, she arrived at the necessity of considering multivalued
(Markov) mappings and suggested me to start a joint research on
multivalued solutions of equations. Our work resulted in a
series of papers, see
\cite{OV1,OV2}. We also had grandiose projects for further research,
for example, of working on
 metric hydrodynamics, but they never were realized.
Approximately at the same time, I started to develop general
(multivalued) dynamics (\cite{V1}), and recently, after a long
interval, I have returned to this subject. In this paper, dedicated to the
unforgettable O.A., I continue this topic.

\section{Markov operators}

\subsection{Definitions}

\begin{definition}
{\it A Markov operator in the Hilbert space
$L^2(X,\mu)$ of complex-valued square-integrable functions on
a Lebesgue--Rokhlin space $(X,\mu)$ with a continuous normalized measure
$\mu$ is a continuous linear operator
$V$ satisfying the following conditions:

{\rm 1)} $V$ is a contraction: $\|V\|\leq 1$ (in the operator norm);

{\rm 2)} $V{\emph{1}}= {V}^*{\emph{1}}={\emph{1}}$, where $\emph{1}$
 is the function identically equal to one;

{\rm 3)} $V$ preserves the nonnegativity of functions:
$Vf$ is nonnegative whenever $f\in L^2(X,\mu)$ is nonnegative}.
\end{definition}

Note that condition 1) follows from  2) and 3), and the second condition in
2) follows from the other ones.

In short: a Markov operator is a
unity-preserving positive contraction.

In the same way we can define a Markov operator
$V$ from one space $L^2(X,\mu)$ to another space $L^2(Y,\nu)$.

\subsection{Three languages}

A geometric analog of a Markov operator
$V:L^2(X,\mu) \to L^2(Y,\nu)$ is a {\it polymorphism}, i.e.,
a measure-preserving multivalued mapping of the space
$(Y,\nu)$ into the space $(X,\mu)$ (for a detailed exposition, see
\cite{V1,V2}).

Each Markov operator is  uniquely $\bmod 0$  generated
by a polymorphism,
and, conversely, for each class of $\bmod 0$ coinciding polymorphisms
there is a Markov operator canonically associated with it.
The correspondence between Markov operators and polymorphisms
extends the classical correspondence
between unitary positive unity-preserving operators and
measure-preserving transformations (Koopman correspondence)
to the case of Markov operators.

This correspondence can easily be explained using the helpful
intermediate notion of bistochastic measure. A bistochastic measure
$\nu$ on the space $X \times X$ is a measure whose projections to the
first and second coordinates (marginal measures) coincide with a given
measure $\mu$. A bistochastic measure is the generalized kernel of a Markov
operator. A bijection between the set of Markov operators
($V$) and the set of bistochastic measures
($\nu$) is given by the formula
$$\nu(A\times B)=\langle V\chi_A,\chi_B\rangle;$$
here $A$ and $B$ are measurable subsets of
$X$ with characteristic functions $\chi_A$ and $\chi_B$, respectively, and
$\nu$ is a bistochastic measure on
$X\times X$, defined on the $\sigma$-field that is the square
of the original $\sigma$-field on
$(X,\mu)$. It is easy to prove that this formula determines a bijection
between Markov operators and bistochastic measures; of course,
$\nu$ can be singular with respect to the product measure
$\mu \times
\mu$, and the conditional measures can be singular with respect to
$\mu$.

In these terms, the polymorphism corresponding to a Markov operator
$V$ is the mapping that associates with
$\mu$-almost every point $x \in X$ some measure on $X$,
namely, the conditional measure
$\mu_V^x=\mu^x $ on the space $X$ regarded as an element of the partition of the space
$(X\times X, \nu)$ into the preimages of points under the projection
from $X\times X$ to the first coordinate;
the conjugate polymorphism associates with a point
$x \in X$ the conditional measure $\mu_x$ on the space $X$
regarded as an element of the partition of the space
$(X\times X, \nu)$ into the preimages of points under the projection
to the second coordinate. By the general Rokhlin theorem,
the conditional measures exist for almost all elements of
a measurable partition.

Each Markov operator
determines a random stationary Markov chain
(discrete-time process)
as follows: its state space
is the original Lebesgue space
$(X,\mu)$, the measure $\mu$ being invariant for this Markov chain,
and the transition probability
$P(\cdot,x)$, $x \in X$, is determined by the kernel of the Markov operator,
more precisely, by the polymorphism that associates with almost every point
$x$ of the space $(X,\mu)$
the conditional measure $\mu^x$ on $X$ (see above); the
transition probabilities are precisely
these conditional measures.
Conversely, each stationary Markov chain determines a bistochastic measure,
namely, the two-dimensional distribution of two adjacent states,
and hence a Markov operator and a polymorphism.

Thus we have three equivalent languages:
\textbf{the language of Markov operators, the language of polymorphisms
and bistochastic measures, and the language of stationary Markov chains}
with continuous state space. For details, see
\cite{V1}. In this paper, we will mainly use the language of Markov operators
and sometimes, when it is helpful, provide explanations in two other languages.

\subsection{Structures}

Markov operators in the Hilbert space
$L^2(X,\mu)$ form a semigroup with respect to multiplication,
with identity element (identity operator), zero element
(one-dimensional orthogonal projection $\theta$ to the subspace of constants),
and involution (operator conjugation $*$).
It is also a convex compact topological semigroup
in the weak
topology in the algebra of all continuous operators.
Indeed, it is easy to verify that the class of Markov operators
is closed under the above-mentioned operations. All these structures
are also defined on the set of polymorphisms of the space
$(X,\mu)$ (or bistochastic measures), and the correspondence
``polymorphism --- Markov operator'' is an antiisomorphism
of semigroups that preserves these structures
(but reverses arrows), see
\cite{V1}.

Since a Lebesgue--Rokhlin space with continuous measure is unique
up to a metric isomorphism (it is isomorphic to
the interval with the Lebesgue measure), the compact space of
Markov operators is also unique in the same sense.
Denote it by
$\cal P$. The subgroup of invertible elements of this semigroup
is precisely the subgroup of positive unitary
unity-preserving operators, i.e., the group of
($\bmod 0$ classes of) all measure-preserving mappings of the Lebesgue
space into itself. Below we give an approximative definition
of the compact space $\cal P$ and all structures on it,
which is independent of operator theory.

\subsection{Approximation lemma}

Consider matrices of order $n$ with nonnegative entries whose
rows and columns sum to $1$ (bistochastic matrices);
they form a convex compact semigroup (with respect to matrix
multiplication), with involution (transposition), zero element
(matrix with all entries equal to
$n^{-1}$), and identity element (identity matrix). Denote by
${\cal P}_n$ the convex compact space of such matrices;
its dimension is equal to $(n-1)^2$;
${\cal P}_n$ is the compact space of Markov operators on
a finite space with the uniform measure.

For positive integers $n$, $m$, $k>1$ with $n=mk$,
partition the rows and columns of matrices from
${\cal P}_n$ into blocks of order $k$. Thus
matrices from ${\cal P}_n$ obtain the structure of block matrices
with blocks of order $k$. Consider the natural projection
$$\pi_{n,m}:{\cal P}_n \rightarrow {\cal P}_m$$
that replaces each matrix block
with the sum of the elements of this block divided by $k$.
The following obvious result, which is however
important for our further considerations, provides another
definition of our main object.

\begin{lemma}
{The compact space $\cal P$ of bistochastic measures,
regarded with all structures defined above (the
structure of a compact space, of a semigroup, etc.),
is the inverse (projective) limit of the spaces
${\cal P}_n$ with respect to the partially ordered set of projections
$(\pi_{n,m}, m|n)$:
$$ {\cal P} ={\underleftarrow\lim}_{n,m}({\cal P}_n,\pi_{n,m}).$$}
\end{lemma}

It is convenient to restrict ourselves only with the powers of a single number
(for instance, $n=2^s$,
$s=1,2,\dots$) and consider the limit along a linearly ordered set.
The proof is straightforward. This construction can easily be interpreted
in terms of weak approximation of operators, but in what follows we will
use this lemma in a slightly different way.

\section{Classes of Markov operators}

\subsection{Generic classes and the list of properties}

The following definition is a specialization of the
well-known terminology.

\begin{definition}
{We say that a class of Markov operators is
\textbf{generic} or forms a set of second category
if this class, regarded as a subset of the set of all Markov operators
$\cal P$, contains an everywhere dense
$G_{\delta}$-set (= intersection of countably many open sets).
A property is called generic if the class of operators satisfying
this property is generic.}
\end{definition}

We will describe generic classes of Markov operators and, consequently,
generic classes of polymorphisms, bistochastic measures, and
stationary Markov chains. As it often happens,
the generic classes we are going to
consider have been scarcely studied, and some
of seemingly paradoxical
generic properties of Markov chains given below
were not even known until recently.

We will be interested in the following properties of Markov operators.

\begin{definition}
A Markov operator $V$ in the space $L^2(X,\mu)$ is called

{\rm 0)} \emph{ergodic} if it has no nonconstant invariant vector; recall that
the spectrum of a Markov operator lies in the unit circle;

{\rm 1)} \emph{mixing} (respectively, \emph{comixing}) if $V^n \to
\theta$ (respectively, ${V^*}^n \to \theta$) as $n \to +\infty$
(recall that $\theta$ is the orthogonal projection to
the subspace of constants);

{\rm 2)} \emph{totally nonisometric} (respectively,
\emph{totally noncoisometric}) if the operator $V$
(respectively, the conjugate operator $V^*$) is
not isometric on any closed invariant subring (sublattice) in
$L^2(X,\mu)$ except that of constants (subrings or sublattices in
$L^2(X,\mu)$ are linear subspaces consisting of functions that are constant
on the elements of some measurable partition of the space
$(X,\mu)$);

{\rm 3)} \emph{dense} if its kernel (= the preimage of zero) is zero and
the image of the space $L^2(X,\mu)$ is a dense linear subspace in
$L^2(X,\mu)$; in other words, if the kernel of $V$ and the kernel
of the conjugate operator
$V^*$ (= cokernel of $V$) are trivial;

{\rm 4)} \emph{extremal} if $V$ (and hence $V^*$)
is an extreme point of the convex compact space $\cal P$ of Markov operators;

{\rm 5)} \emph{indecomposable} if there is no measurable subset
$A \subset X$ of positive measure $\mu$ and a measurable subset
$B \subset X\times X$ of positive measure $\nu$ such that
the image of the characteristic function $\chi_B$ under the operator
$V$ and the conjugate operator
$V^*$ is positive and strictly less than $1$ almost everywhere on the set
$A$.
\end{definition}

The zero, first, and second properties should be called
operator-dynamical; the third one, measure-geometric;
the fourth and the fifth ones, properly geometric.

\subsection{Analysis of properties}

Let us comment on the above definitions.

\smallskip
1. The term and notion of mixing appeared in the theory
of dynamical systems; here it means that the shift in the space
of trajectories of the corresponding stationary Markov chain
is a mixing in the sense of that theory.
In the case of Markov chains {\it with finite state space},
this property is equivalent
to a much stronger property of a chain --- the triviality of
the tail $\sigma$-field at minus infinity,
which means that
the $\sigma$-field $\bigcap_{n=0}^{\infty} {\cal
A}_{-\infty}^{-n}$ of measurable subsets in the space of two-sided
trajectories of a Markov chain consists of two elements:
the class of zero-measure sets and the class of the whole space.
This property has many other names and many other formulations
(Kolmogorov regularity, 0--2 law, etc.)
\cite{R}. In some cases, there are well-known conditions of mixing;
for example, an aperiodic chain with finitely or countably many
states and, more generally, an aperiodic chain satisfying the Harris condition
(see \cite{R}) are mixing. The notions of mixing and comixing
are in general position.

In the general theory of contractions in Hilbert spaces
(see \cite{FN}), one considers four classes of contractions,
depending on whether or not
the sequence of positive powers
of the operator or the conjugate operator weakly tends to zero.
Depending on what of the four variants takes place, one uses
the notation
$C_{0,0}$, $C_{0,1}$, $C_{1,0}$, $C_{1,1}$. Borrowing this notation,
we can say that the class
of nonmixing and noncomixing Markov operators,
which is most important for our purposes, lies in
$C_{1,1}$.

\smallskip
2. The question arises with the notions of mixing and totally nonisometry: whether every totally
nonisometric Markov operator is mixing. One of the main points of the general theory of
contractions (see \cite{FN}) is that for contractions it is not the case; in other words, there
exists a totally nonisometric (respectively, noncoisometric) contraction such that the sequence of
positive (respectively, negative) powers does not tend to zero. It turns out that in the theory of
Markov operators such an effect also takes place: there exist nonmixing and noncomixing Markov
operators that are still totally nonisometric and noncoisometric. In order to understand the
paradoxical nature of this situation, let us reformulate the condition of being totally
nonisometric in geometric terms.

\begin{definition}
{A stationary Markov chain with state space $X$ and invariant measure $\mu$
is called {\it totally nondeterministic} if there is no
measurable partition
$\xi$ of the space $(X, \mu)$ such that the transition operator
acts deterministically
on its blocks, i.e., sends a block to a block.}
\end{definition}

\begin{lemma}
{A Markov chain is totally nondeterministic if and only if the corresponding
Markov operator is totally nonisometric.}
\end{lemma}

In the theory of chains with finitely many states, the property
of being totally nondeterministic is known as
``the absence of subclasses'' or ``aperiodicity,'' etc., see
\cite{R}.  In
\cite{V1}, for certain reasons,
the corresponding polymorphisms were
called {\it simple}. For chains with countably many states and,
more generally, for chains satisfying the Harris condition
(see \cite{R}), the condition of being aperiodic, i.e., totally nonisometric,
is equivalent to mixing.

However, for general Markov chains this is not the case: there exist
totally nonisometric and nonmixing Markov operators, i.e., totally
nondeterministic and nonmixing Markov processes. The first
example of this type is due to M.~Rosenblatt \cite{Ro}. For other
examples, see \cite{V1,V3}; in these examples, the transition probabilities
are singular with respect to the invariant measure and their
behaviour is rather complicated; the behaviour of the powers of the Markov
operator is also quite nontrivial. Here we do not describe
these examples, referring the reader to the above-mentioned papers,
but we will prove that they are generic.

\smallskip
3. Property 3) need no comments; it means that the operator
may have no bounded inverse, yet the inverse operator
exists on an everywhere dense subset, and the same is true for
the conjugate operator. However, it is worth explaining why
this property is important.
Let us say that a Markov operator
$V$ is a {\it quasi-image} (in \cite{FN}, it was called
a {\it quasi-affinitet})
of a Markov operator $W$ if there exists a dense Markov operator
$L$ such that $$LV=WL.$$
If $V$ is a quasi-image of
$W$ and $W$ is a quasi-image of
$V$ (i.e., there exists a dense Markov operator
$M$ such that $MW=VM$), then we will say that the operators
$V$ and $W$ are
{\it quasi-similar}. Quasi-similarity is an equivalence relation;
it would not be of any interest if we did not require that the intertwining
operators $L$ and $M$ should be dense. In the theory of contractions
there is a number of important results on quasi-similarity
(see \cite{FN}), but it seems that for Markov operators this notion
has never been introduced and studied. The main problem, which we do
not discuss here, is {\it when two unitary Markov operators
are quasi-similar and what unitary Markov operators can be quasi-similar
to totally nonisometric Markov operators}. These problems are
extremely important for the theory of dynamical systems and
statistical physics in connection with discussion on
irreversibility (see references in
\cite{V2, V3}).

\smallskip
4, 5. The notions of extremality and indecomposability are of
completely different nature. If a Markov operator is a
nontrivial convex combination of other Markov operators, this means
that the corresponding Markov shift is a
{\it skew product over a Bernoulli shift}; in other words, it is
a random walk over the trajectories of Markov components. In particular,
if the Markov operators occurring in the convex combination are unitary,
then we have a
{\it random walk over the trajectories of deterministic
transformations with invariant measure} or a so-called
random dynamical system. This is just the case for Markov chains
with finitely many states and the uniform measure, because, by the
Birkhoff--von Neumann theorem, extreme points of the polyhedron
of bistochastic matrices are permutation matrices.
In the case of general bistochastic measures, extremal Markov operators
are not necessarily unitary; moreover, the conditional measures can
even be continuous (see below). From the geometric
point of view, bistochastic measures  were studied by many authors; see,
e.g., \cite{S,V2} for nontrivial examples of extremal
polymorphisms and Markov operators and further references.
An obvious necessary condition for extremality is as follows: {\it there is no
set of constant width strictly between $0$ and $1$ with respect to both
projections}, i.e., there is no measurable set of
intermediate measure such that the images of the characteristic function
of this set under the Markov operator and its conjugate are constant
functions. We may go further and introduce the notion
of indecomposability (see above). If an operator is indecomposable,
then it cannot be represented as a convex
combination of other Markov operators, even with nonconstant
(depending on the point) coefficients that are not equal to $0$ or $1$
at sets of positive measure. In this case, the shift in
the space of trajectories of the corresponding Markov process
cannot be represented as a random walk over the trajectories
of any Markov shifts with probabilities depending on
the point and different from $0$ and $1$ almost everywhere. It turns out that
even this condition, which is much stronger than the usual extremality,
determines a generic class of Markov operators. The indecomposability means
that there is {\it no subset of nontrivial width with respect to both projections
over the whole space or at least over a set of positive measure}.
It is not difficult to deduce that indecomposability implies
extremality, but the converse is not true. A remarkable characteristic
property of every indecomposable bistochastic measure $\nu$
is that in the space
$L^2(X\times X, \nu)$ every function can be approximated by functions of
the form $f(x)+f(y)$; in other words, {\it there is no nonzero function that has
zero expectation with respect to both subalgebras}.

\section{The main theorem}

\subsection{Formulation}
We consider a Lebesgue--Rokhlin space
$(X,\mu)$ with continuous measure. All Markov operators act in the space
$L^2(X,\mu)$ of square-integrable complex-valued functions. Recall that
a Markov operator is ergodic if it has no
nonconstant invariant vector. Below, by singularity we mean the
singularity with respect to the measure $\mu$.

\begin{theorem}
{A generic Markov operator enjoys the following properties:

{\rm 1)} its spectrum has no discrete component (in the orthogonal
complement to the subspace of constants); in particular, it is ergodic;

{\rm 2)} it is neither mixing nor comixing;

{\rm 3)} it is totally nonisometric and totally noncoisometric;

{\rm 4)} it is dense;

{\rm 5)} it is extremal and indecomposable;

{\rm 6)} almost all its transition probabilities are continuous
and singular.}
\end{theorem}

\smallskip\noindent{\it Remark.}
Most papers on the theory of Markov chains deal with the cases
of either absolutely continuous or discrete transition
probabilities (for example, Doeblin condition, Harris condition, etc.).
In these cases, it is difficult to discover
most important and generic effects, such as the absence of mixing
for totally nondeterministic operators, as well as other generic properties.
\smallskip

\subsection{Proof}

Since the intersection of finitely or countably many generic
classes in a complete metrizable separable space is generic
(Baire theorem), it suffices to prove that each class is generic by itself.
Further, all ``coproperties'' are similar
to the corresponding properties for the conjugate operator, and hence
they are generic provided that the original properties are generic.

We start by observing
that the group of positive unitary operators
is everywhere dense in $\cal P$; this follows from the approximation lemma
and the following simple fact: every rational bistochastic matrix
of order $n$ whose all entries have
denominator $N$ is the projection
($\pi_{Nn,n}$) of some permutation matrix of order
$nN$ (see \cite{V2}). This implies that the sets of operators
satisfying properties
1, 2, 4, 5 are everywhere dense, because a generic measure-preserving automorphism
is ergodic, has a simple continuous spectrum (\cite{Yuz}), is nonmixing
(since it is deterministic), extremal, and, of course, indecomposable
and dense. The fact that the set of totally nonisometric operators is
everywhere dense and even satisfies the Baire property
($G_{\delta}$) also follows from the lemma,
but in this case we should use other  matrices,
namely, irreducible ones: irreducible bistochastic matrices
are generic even in the finite-dimensional case, and the
projections preserve irreducibility.

Let us verify that the sets of operators satisfying the
remaining properties are  $G_{\delta}$-sets.
The fact that property 6 is generic can also be seen from the lemma, because
the existence of a discrete component in the conditional measures can
be written in terms of the approximating bistochastic matrices
(see \cite{V2} for details). Extreme points of every convex compact set
form a $G_{\delta}$-set (see \cite{F}). The $G_{\delta}$-condition
for indecomposable bistochastic measures is proved in
\cite{V2}.  For properties 2 and 4, the
$G_{\delta}$-condition is trivial. Finally, this condition
for property 1
is satisfied in the algebra of all bounded operators: the set of
operators that have no
discrete component in the spectrum is a
$G_{\delta}$-set; and so is its intersection with the set of Markov operators.

\subsection{Remarks, problems, conjectures}

1. As it often happens,
it is easier to prove that a property is generic than to construct
an explicit example of a generic object.
The deepest problem, which has
a nontrivial solution, the construction of a Markov
operator that simultaneously satisfies properties 2 and 3, i.e., is
nonmixing and totally nonisometric, is considered in the recent author's
paper \cite{V1}. In that paper, a relation is established of such examples to the
hyperbolic theory of dynamical systems. This led to a new
characterization of $K$-systems and to the following conjecture,
which we state here without specifying details (this would require
giving new definitions and  will be done elsewhere):

\smallskip
{\bf A generic polymorphism is a singular random perturbation of
a Kolmogorov automorphism. Correspondingly, a generic Markov operator
is a singular perturbation of a unitary positive operator
conjugate to a $K$-automorphism.}
\smallskip

It seems possible to construct a Markov operator (polymorphism)
that simultaneously satisfies all properties 1--6 by specializing
these new examples.

\smallskip
2. It is natural to ask whether the class of Markov operators such that
the shift in the space of trajectories
of the corresponding Markov chain is a $K$-automorphism is generic. We emphasize that,
as observed above, for a generic Markov operator, the Markov
generator is not a $K$-generator, because there is
no mixing. However, in all known examples, there exists another,
non-Markov, $K$-generator.

The next question: what Markov operators generate a shift that
is (isomorphic to) a Bernoulli shift? Unfortunately, Ornstein's technique
($\bar d$-metric) is not suitable for studying
processes with
continual state space. The famous Kalikow's example
\cite{Kal} of a non-Bernoulli and even non-loosely Bernoulli automorphism
(a random walk over $(T,T^{-1})$, where $T$ is a Bernoulli shift)
demonstrates the wide possibilities
of the natural Markov generators.
A close and, apparently, difficult question is whether every
$K$-automorphism has a Markov $K$-generator.

\smallskip
3. There is an acute problem concerning the definition of the entropy of a Markov
operator or polymorphism. There are different suggestions,
and, probably, there should be ``different'' entropies
corresponding to different properties of polymorphisms.
One definition is introduced, not quite distinctly, in
\cite{V4}; this entropy is positive for generic polymorphisms of a
finite space (bistochastic matrices). Another definition
of the entropy of Markov operators, as well as further
references, can be found in \cite{DF}.
In ergodic theory, the generic
value of the Kolmogorov entropy is zero. The answer to the same question
for polymorphisms and any of the entropies is open.

\smallskip
4. Every polymorphism (Markov operator) generates certain equivalence
relations on $(X, \mu)$. One of them is the partition into orbits:
two points $x$ and $y$ lie in the same orbit if there exist positive
integers $n,m$ such that the conditional measures
$\mu_{V^n}^x$ and $\mu_{V^m}^y$ are not singular (their mutual densities
are not identically zero or infinite). Another one is the
transitive envelope
of the nonsingularity relation for conditional measures:
$x \sim y$ if there exists a positive integer $k$ and a finite chain of points
$x_0=x, x_1, \dots, x_k=y$ such that the conditional measures
$\mu_V^{x_i}$ and $\mu_V^{x_{i+1}}$ are not singular for
$i=0, \dots, k-1$. What is the generic
behaviour of these equivalence relations?
It seems that the class of Markov operators for which these equivalence relations
are ergodic is generic.

\smallskip
5. We will not discuss other generic properties of Markov operators.
Let us just mention one important link to the theory of groupoids,
and one more problem.
Consider a Markov operator $V$, the conjugate operator
$V^*$, and the semigroup spanned by these two operators. Let
$T=\sum_w{c_w w(V,V^*)}$, where $w(\cdot,\cdot)$ ranges over all
finite words in the alphabet $\{V,V^*\}$ and $c_w$
are positive coefficients decreasing fast enough
with total
sum equal to $1$. The operator $T$ is again a Markov operator; it
is related to the measurable groupoid generated by the original operator,
more precisely, by the first equivalence relation mentioned above, the
partition of the space $(X, \mu)$ into orbits. It is very important to find out
when this partition is hyperfinite (in another terminology,
tame) and whether this case is generic; most likely, it is not.

\smallskip
6. On the other hand, with every Markov operator we can associate
the $C^*$-algebra generated by this operator, the conjugate operator, and
the operators of multiplication by some class of bounded
measurable functions. Such algebras generalize the notion of skew
product, and the study of their properties (for example, amenability,
simplicity, etc.) will give new examples of $C^*$-algebras.

\smallskip
7. Finally, the last question is also related to the theory of
$C^*$-algebras. Of great interest is the
$C^*$-algebra generated by all Markov operators. It seems that is has never
been considered. This algebra is not separable and does not coincide
with the algebra of all operators. Presumably, it coincides with
the algebra of all operators that preserve order-bounded
sets in $L^2$, and its elements have a natural integral representation
with a kernel that may be nonpositive.
It must play the same role in the theory of Markov operators
and dynamical systems that the algebra of all operators plays
in general operator theory.

\smallskip
8. One particular class of polymophisms is of special interest, namely, the
class of so-called {\it algebraic polymorphisms},
i.e., correspondences on compact Abelian
groups with the Haar measure. A typical example is an
algebraic polymorphism of the circle, i.e., the uniform
measure on the one-dimensional cycle of the two-dimensional
torus determined by the equation
$u^p=v^q$, where  $p,q$ are positive integers and
$u,v$ are coordinates on the torus. These examples will
be considered from different viewpoints in a joint paper
by the author and K.~Schmidt which is now in preparation.

\section{Comments}

\subsection{Relation to previous works}

In the joint papers with O.~A.~Ladyzhenskaya
\cite{OV1,OV2} mentioned in the introduction, we considered
application of theorems on measurable selection. In order to
apply general theorems,
we needed some {\it a priori} estimates. More or less
simultaneously, there appeared several papers in the same direction ---
Foias, Temam, and others --- with similar results.
The common feature of all these works was that they
regarded multivaluedness
out of the context of Markov operators, i.e.,
with Cauchy data one associated the set of solutions rather than a
measure on them. The theory of polymorphisms and Markov operators
can have more subtle applications and simulate more complicated phenomena
than theorems on measurable selection.
However, for that one needs
to develop the theory of one-parameter semigroups of Markov operators
and polymorphisms and their Lie generators. Apparently,
the above results on generic properties can be generalized
to the case of such semigroups, but this is still to be done.
One can hope that generic semigroups of polymorphisms
will also find applications in hydrodynamics, as we once
discussed with O.A. I have been keeping the memory of
our discussions through all these years.

\subsection{Acknowledgments}

The author repeatedly
discussed the theory of Markov operators with
M.~I.~Gordin, who owns important ideas
relating hyperbolic
theory to probability theory; he also gave me
several useful references.
Many years ago, D.~Z.~Arov drew my attention to
the tentative relation between
the theory of polymorphisms and the Foias--Nagy
theory of contractions, and at the same time, B.~Rubshtein
pointed out M.~Rosenblatt's example to me.
With V.~A.~Kaimanovich, I discussed
the relation to the theory of boundaries,
and with L.~A.~Khalfin
(1932--1998), we talked about possible links to physics. To all of them
I am deeply grateful.

\end{document}